\documentclass[10pt]{article}
\usepackage{amsfonts}
\usepackage{amsmath}
\usepackage{amssymb}
\usepackage[margin=1in]{geometry}

\def\R{\mathbb{R}}
\def\N{\mathbb{N}}

\def\Z{\mathbb{Z}}

\def\nl{\newline}
\def\bo{\nl\phantom{a}\hfill $\Box$\nl}

\makeatletter
\@namedef{subjclassname@2020}{\textup{2020} Mathematics Subject Classification}
\makeatother

\numberwithin{equation}{section}

\renewcommand{\Re}{\operatorname{Re}}

\newcommand{\eps}{\varepsilon}

\newcommand{\MP}{\mathcal{P}}
\newcommand{\MN}{\mathcal{N}}
\newcommand{\MF}{\mathcal{F}}

\numberwithin{equation}{section}

\begin{document}
\begin{center}
{\Large {\bf A Mean value Theorem for general Dirichlet series}}\vspace{0.1in}\\
 Frederik Broucke and Titus Hilberdink
\end{center}
\indent
\begin{abstract} In this paper we obtain a mean value theorem for a general Dirichlet series $f(s)= \sum_{j=1}^\infty a_j n_j^{-s}$ with positive coefficients for which the counting function $A(x) = \sum_{n_{j}\le x}a_{j}$ satisfies $A(x)=\rho x + O(x^\beta)$ for some $\rho>0$ and $\beta<1$. We prove that $\frac1T\int_0^T |f(\sigma+it)|^2\, dt \to \sum_{j=1}^\infty a_j^2n_j^{-2\sigma}$ for $\sigma>\frac{1+\beta}{2}$ and obtain an upper bound for this moment for $\beta<\sigma\le \frac{1+\beta}{2}$. We provide a number of examples indicating the sharpness of our results.\nl

\vspace{0.3in}\end{abstract}

\noindent
{\large {\bf 1. Introduction}}\nl
Given a function defined by a Dirichlet series, $f(s)=\sum_{n=1}^\infty a_n n^{-s}$, the existence of the mean value, namely
\[ \lim_{T\to\infty} \frac1{2T}\int_{-T}^T|f(\sigma+it)|^2\, dt,\]
is a well-known object of study and has received much attention. First and foremost this was motivated by attempts to understand the Riemann zeta function $\zeta(s)$ and more generally $L$-functions in the critical strip. The mean value and indeed higher moments of $\zeta(s)$ play an important role in analytic number theory. The study of these moments has a rich history dating back to the beginning of the twentieth century with notable contributions by for instance Landau \cite[Teil 24]{Landau}, Hardy and Littlewood \cite{HardyLittlewood}, and Ingham \cite{Ingham}. However, our understanding of these moments is still very much incomplete, and this remains an area of ongoing research (see e.g.\ \cite{Ivic, Sound, Harper} and the references therein).

More generally, series of the form $\sum_{j=1}^{\infty}a_{j}n_{j}^{-s}$, called \emph{general Dirichlet series}, have also received much attention (see e.g.\ \cite{Landau, HardyRiesz, CI, CII, Potter}). Here $(n_{j})_{j}$ is a positive strictly increasing unbounded sequence.  
Our original motivation for this article has come from studying the \emph{Beurling zeta function} $\zeta_{\MP}(s)$ associated to a Beurling prime system. Here
\[	
	\zeta_{\MP}(s) = \sum_{j=1}^{\infty}\frac{a_{j}}{n_{j}^{s}} = \prod_{j=1}^{\infty}\frac{1}{1-p_{j}^{-s}},
\]
where $(\MP, \MN)$, $\MP = (p_{j})_{j\ge1}$, $\MN = (n_{j})_{j\ge1}$ is a system of \emph{Beurling generalized primes and integers} \cite{B, DiamondZhang}, with $a_{j}$ denoting the multiplicity of $n_{j}$.  Recently, the second moment of Beurling zeta functions has been studied in \cite{DGN}. There, the existence of the mean value was established in certain cases, but under the assumption of a certain separation condition on the Beurling integers $n_j$. In this article, we study the mean value of a general class of Dirichlet series, without assuming any separation of the $n_j$.

Our setup is the following: let $(n_{j})_{j\ge1}$ be a strictly increasing unbounded sequence of positive real numbers, and let $a_{j}>0$ for $j\ge1$. Set $A(x) = \sum_{n_{j}\le x}a_{j}$, and suppose that 
\[	A(x) = \rho x + O(x^{\beta}),\tag{1.1}\]
for some $\rho>0$ and $0\le\beta<1$. The corresponding Dirichlet series, 
\[
	f(s) = \sum_{j=1}^{\infty}\frac{a_{j}}{n_{j}^{s}},
\]
converges absolutely for $\Re s = \sigma>1$.  In view of (1.1), $f(s)-\rho/(s-1)$ admits analytic continuation to the half-plane $\sigma>\beta$. We are interested in the asymptotic behavior of 
\[
	\frac{1}{T}\int_{0}^{T}|f(\sigma+it)|^{2}\, dt \tag{1.2}
\]
as $T\to\infty$. 
By absolute convergence of the Dirichlet series the limit exists if $\sigma>1$ and is given by $\sum_{j=1}^{\infty}\frac{a_{j}^{2}}{n_{j}^{2\sigma}}$. The question is, how far to the left does the mean value exist? In this article, we show that the mean value exists if $\sigma>\frac{1+\beta}{2}$, and that this is sharp. Second, we also show a sharp upper bound for (1.2) in the strip $\beta<\sigma\le \frac{1+\beta}{2}$. Hence, we provide a complete solution to the problem of determining the maximal growth of the mean square of such general Dirichlet series.\nl

\noindent
{\bf Theorem 1.1}\nl
{\em Let $(n_{j})_{j\ge1}$, $(a_{j})_{j\ge1}$ be as above and suppose that $(1.1)$ holds.
\begin{enumerate}
	\item For\footnote{For $\sigma=1$, we start integrating at $t=1$ to avoid the pole.} $\sigma>\frac{1+\beta}{2}$  
\[ \lim_{T\to\infty}\frac1T\int_0^T |f(\sigma+it)|^2 \, dt = \sum_{j=1}^\infty\frac{a_j^2}{n_j^{2\sigma}}.	\]
	\item  Uniformly for $\beta+\delta\le\sigma\le\frac{1+\beta}{2}$ (any $\delta>0$), we have
\[
	\frac{1}{T}\int_{0}^{T}|f(\sigma+it)|^{2}d t \ll_\delta \frac{T^{\frac{1+\beta-2\sigma}{1-\beta}}-1}{1+\beta-2\sigma},
\]
and where the right hand side is taken to be $\log T$ when $\sigma=\frac{1+\beta}{2}$.
\end{enumerate} }
\noindent
{\bf Remark 1}\, (i)\,  It follows immediately that the mean value exists for $\sigma>\frac{1+\beta}{2}$ under the slightly weaker $A(x) =\rho x+O_{\eps}(x^{\beta+\eps})$ for every $\eps>0$. 

(ii)\, Note that $a_{j} \ll n_{j}^{\beta}$, so that 
\[
	\sum_{j=1}^{\infty}\frac{a_{j}^{2}}{n_{j}^{2\sigma}} \ll \sum_{j=1}^{\infty}\frac{a_{j}}{n_{j}^{2\sigma-\beta}},
\]
which converges for $\sigma>\frac{1+\beta}{2}$. Bessel's inequality, adapted to general Dirichlet series (see the appendix), gives
\[ 
	\liminf_{T\to\infty}\frac{1}{T}\int_0^T |f(\sigma+it)|^2\, dt \ge \sum_{j=1}^{\infty}\frac{a_{j}^{2}}{n_{j}^{2\sigma}}. \tag{1.3}
\]
 Thus the convergence of the series on the right is a necessary condition for a mean value to exist. As the RHS above may diverge for $\sigma \le \frac{1+\beta}{2}$, we cannot, in general, have a mean value for such $\sigma$, as we shall see. We shall also see that convergence of the series in (1.3) is not sufficient for $f$ to have a mean value.

(iii)\, Applying Theorem 1.1 to the Riemann zeta function, we obtain that $|\zeta(\sigma+it)|^2$ has a mean for any $\sigma>\frac12$, simply from the fact that $[x]=x+O(1)$  (or even $O(x^\eps)$ for any $\eps>0$). This is in contrast to all the known proofs of this result which use some further knowledge of the spacing of the integers, either directly as in Landau \cite[p. 806]{Landau} (see also Titchmarsh \cite[Lemma 7.2]{Titchmarshzeta}\footnote{He shows $\sum_{m<n<T} \frac1{(mn)^\sigma\log\frac{n}{m}}\ll T^{2-2\sigma}\log T$, where $m,n$ range over $\N$. If, for example, $n-m=e^{-n}\sim e^{-T}$ were possible, then the sum would be exponentially large.}), or indirectly such as via van der Corput's method for exponential sums, the Fourier series of $[x]$, or some general form of Hilbert's inequality (see also the Montgomery--Vaughan Theorem below).

(iv)\, 
In \cite[Theorem 1]{H}, it was proven that if the $n_j$ are the integers of a Beurling prime system\footnote{One can show this result extends to our setting with the condition that $n_j\gg j$.} with multiplicities $a_{j}$ satisfying (1.1) with 
$\beta<\frac12$, then for $\sigma\in (\beta,\frac12)$, (with $f(s)=\zeta_\mathcal{P}(s)$)
\[ \int_0^T |f(\sigma+it)|^2\, dt  =\Omega_{\eps}(T^{2-2\sigma-\eps})\]
for all $\eps>0$.
When $\beta=0$, Theorem 1.1(b) says that 
\[ 
	\int_{0}^{T}|f(\sigma+it)|^{2}d t = O(T^{2-2\sigma})\quad \mbox{ for $\sigma\in (0,1/2)$, and } \quad \int_{0}^{T}\Bigl|f\Bigl(\frac{1}{2}+it\Bigr)\Bigr|^{2}d t = O(T\log T).
\]
This shows that for $\beta=0$, the mean square is quite rigid: for $\sigma > 1/2$, it is asymptotic to $c_{\sigma}T$, for $\sigma=1/2$ its growth is between $T$ and $T\log T$, and for $0<\sigma<1/2$, its growth is between $T^{2-2\sigma-\eps}$ (every $\eps>0$) and $T^{2-2\sigma}$.

Another natural question is whether 
\[ 
	\lim_{T\to\infty} \frac1{T^{2-2\sigma}}\int_{0}^{T}|f(\sigma+it)|^{2}d t\quad\mbox{ and }\quad \lim_{T\to\infty}\frac1{T\log T}\int_{0}^{T}\Bigl|f\Bigl(\frac{1}{2}+it\Bigr)\Bigr|^{2}d t
\]
exist. We show by example in Section 4 that they need not exist, even for ordinary Dirichlet series.\nl

As mentioned before, the topic of general Dirichlet series was a popular area of study in the first half of the twentieth century (see e.g.\ \cite{Landau, HardyRiesz, CI, CII, Potter}). Regarding mean values, it is worth mentioning a result of Landau \cite[Satz 37, p.\ 793]{Landau}, which is somewhat related to our Theorem 1.1. This result requires a certain separation of the sequence $(n_{j})_{j}$ known as Landau's condition:
\[
	\forall \delta > 0: n_{j+1} - n_{j} \gg_{\delta} e^{-n_{j}^{\delta}} \tag{LC}.
\]
Landau showed: \nl

\noindent
{\bf Theorem} [Landau]\nl
{\em Suppose the sequence $(n_{j})_{j\ge1}$ satisfies (LC). If the Dirichlet series $g(s) = \sum_{j=1}^{\infty}c_{j}n_{j}^{-s}$ with complex coefficients has abscissa of convergence $\sigma_{c}$ and abscissa of absolute convergence $\sigma_{a}$, then
\[
	\lim_{T\to\infty}\frac{1}{2T}\int_{-T}^{T}|g(\sigma+it)|^{2}dt = \sum_{j=1}^{\infty}\frac{|c_{j}|^{2}}{n_{j}^{2\sigma}}, \quad \text{for } \sigma > \frac{\sigma_{c}+\sigma_{a}}{2}.
\]
}
Separation conditions such as (LC) are in fact quite common in the study of general Dirichlet series. We emphasize that no separation is required for our Theorem 1.1.
As our proofs do not require modern techniques, it may be a little surprising that a result like ours was not obtained previously. 

\medskip
The paper is organized as follows. In Section 2, we prove Theorem 1.1 by means of three propositions. Next in Section 3, we provide several examples showing the sharpness of Theorem 1.1. In Section 4, we investigate the behavior of $(1.2)$ on the ``critical line'' $\sigma = \frac{1+\beta}{2}$. Assuming additionally a certain separation property of the $n_{j}$, we can employ a result of Montgomery--Vaughan to show an asymptotic formula (Theorem 4.1). 
Finally in Section 5, we apply Theorem 1.1 to bound the number of zeros of $f(s)$ to the right of $\sigma=\frac{1+\beta}{2}$. 

We also include an appendix, where we briefly discuss an extension of Theorem 1.1 to general non-decreasing right-continuous functions $A(x)$, and where we provide some details on $(1.3)$.

\medskip

\noindent
{\bf Notation}: The various symbols $O, o, \ll, \Omega, \sim$ have their usual meaning. For real functions $f,g$ defined on neighbourhoods of infinity, $f=O(g)$ (equivalently, $f\ll g$) means $|f(x)|\le Cg(x)$ for some constant $C$ and all $x$ sufficiently large. We write $f=O_\eps(g)$ (or $f\ll_\eps g$) to indicate the constant may depend on a parameter $\eps$; $f=o(g)$ means $f(x)/g(x)\to 0$ as $x\to \infty$, $f\sim g$ means $f(x)/g(x)\to 1$ as $x\to \infty$, while $f=\Omega(g)$ if $f\neq o(g)$. Also we use $f\prec g$ (equivalently $g\succ f$) to mean $f=o(g)$.
\bigskip

\noindent
{\large {\bf 2. Proofs}}\nl
Before proving Theorem 1.1 in full generality, we mention how Landau's theorem and knowledge of the mean value of the Riemann or Hurwitz zeta functions can yield Theorem 1.1(a) under (LC). This was in fact the way in which the authors of \cite{DGN} showed existence of the mean value of certain Beurling zeta functions. We refer to their paper for more details.

Given a Dirichlet series with positive coefficients $f(s) = \sum_{j=1}^{\infty}a_{j}n_{j}^{-s}$ satisfying (1.1) and (LC), one considers an auxiliary Dirichlet series $h(s) = \sum_{j=1}^{\infty}m_{j}^{-s}$, where the $m_{j}$ are such that $\sum_{m_{j}\le x}1 =  x + O(1)$ and $(n_{j}) \cup (m_{j})$ still satisfies (LC), and such that we know the existence of the mean value of $h(s)$ for $\sigma>1/2$. For example, one may take $m_{j} = j + \alpha_{j}$ for a suitable sequence $\alpha_{j}\in [0,1)$, and compare to a Hurwitz zeta function to establish the mean value of $h(s)$. Then $g(s) = f(s) - \rho h(s)$ satisfies the requirements of Landau's theorem with $\sigma_{c} = \beta$ and $\sigma_{a} = 1$, so we obtain 
\[
	\frac{1}{T}\int_{0}^{T}|f(\sigma+it)|^{2}dt \ll 1 \quad \text{for } \sigma > \frac{1+\beta}{2},
\]
which is enough to conclude $\lim_{T\to\infty}\frac{1}{T}\int_{0}^{T}|f|^{2} = \sum_{j=1}^{\infty}a_{j}^{2}n_{j}^{-2\sigma}$ by a classical result of Carlson \cite{CI}. 

\bigskip
We now return to our general setup (i.e.\ without a separation condition), and give the proof of Theorem 1.1. We write 
\[ f_N(s) = \sum_{n_j\le N} \frac{a_j}{n_j^s}\]
for the partial sums. 
\nl

\noindent
{\bf Proposition 2.1}\nl
{\em Suppose $(1.1)$ holds. Then, for $\sigma\in (\frac{1+\beta}{2},1]$, 
\[ \lim_{T\to\infty}\frac{1}{T}\int_0^T |f_N(\sigma+it)|^2\, dt = \sum_{j=1}^\infty\frac{a_j^2}{n_j^{2\sigma}}\]
as $N,T\to\infty$ such that $T\succ N^{2-2\sigma}\log N$ if $\sigma<1$ and $T\succ(\log N)^2$ if $\sigma=1$.}\nl

\noindent
{\em Proof.}\, Fix $\sigma\in (\frac{1+\beta}{2},1)$. We deal with $\sigma=1$ after. We have
\begin{align*}
\frac1T \int_0^T |f_N(\sigma+it)|^2\, dt &=\frac1{2T} \int_{-T}^T \biggl|\sum_{n_j\le N}\frac{a_j}{n_j^{\sigma+it}}\biggr|^2\, dt = \sum_{n_j,n_k\le N}\frac{a_ja_k}{(n_jn_k)^\sigma}\frac1{2T} \int_{-T}^T \Bigl(\frac{n_j}{n_k}\Bigr)^{it}\, dt \\
& =\sum_{n_j\leq N}\frac{a_j^2}{n_j^{2\sigma}} + \sum_{n_j<n_k\le N}\frac{a_ja_k}{(n_jn_k)^\sigma}\frac1{2T} \int_{-T}^T \Bigl(\frac{n_j}{n_k}\Bigr)^{it}+\Bigl(\frac{n_k}{n_j}\Bigr)^{it}\, dt\\
& =\sum_{n_j\leq N}\frac{a_j^2}{n_j^{2\sigma}} + 2\sum_{n_j<n_k\le N} \frac{a_ja_k s(T\log\frac{n_k}{n_j})}{(n_jn_k)^{\sigma}},\tag{2.1}
\end{align*}
where $s(x)=\frac{\sin x}{x}$ (with $s(0)=1$). As $N\to\infty$, the first term converges to the desired limit, so we must show the double sum converges to 0 if $T\succ N^{2-2\sigma}\log N$. We use the elementary inequality 
\[ |s(x)|\le \min\Bigl\{1,\frac1x\Bigr\}.\]
First note that for a given $n_0$, the terms with $n_k\le n_0$ in (2.1) contribute at most $O_{n_{0}}(\frac1T)$. Thus we need only consider $n_k\ge n_0$ in (2.1). 

For $n_j>n_k-n_k^\beta$, we use $|s(x)|\le 1$. This part of the double sum in (2.1) is, in modulus, at most
\[ 2\sum_{n_0\le n_k\leq N} \frac{a_k}{n_k^{\sigma}}\sum_{n_k-n_k^\beta\le n_j<n_k}\frac{a_j}{n_j^{\sigma}}\ll \sum_{n_0\le n_k\leq N} \frac{a_k(A(n_k)-A(n_k-n_k^\beta))}{n_k^{2\sigma}}\ll \sum_{n_k\ge n_0} \frac{a_kn_k^\beta}{n_k^{2\sigma}}.\]
As $2\sigma-\beta>1$, the series converges and, by taking $n_0$ large enough, this part of the sum can be made as small as we please.

For $n_j\le n_k-n_k^\beta$, we use $|s(x)|\le \frac1x$ and the inequality $\log(\frac{n_k}{n_j})\ge \frac{n_k-n_j}{n_k}$. This part of the sum is, in modulus, at most
\[\frac2T\sum_{n_k\leq N} \frac{a_k}{n_k^{\sigma}}\sum_{n_j\le n_k-n_k^\beta} \frac{a_jn_k}{n_j^{\sigma}(n_k-n_j)} \ll \frac1T\sum_{n_k\leq N} \frac{a_k}{n_k^{\sigma}}\biggl(\sum_{n_j\le \frac{n_k}{2}}\frac{a_j}{n_j^\sigma}+n_k^{1-\sigma}\sum_{\frac{n_{k}}{2}<n_j\le n_k-n_k^\beta} \frac{a_j}{n_k-n_j}\biggr).\tag{2.2}\]
The part with $n_j\le\frac{n_k}{2}$ is\footnote{Here we use the fact that $\sum_{n_j\leq x} a_jn_j^{\lambda} \sim \frac{\rho}{1+\lambda}x^{1+\lambda}$ for fixed $\lambda>-1$.} $\ll \frac1T\sum_{n_k\le N} a_kn_k^{1-2\sigma}\ll \frac{N^{2-2\sigma}}{T}\to 0$. 

We must estimate the inner sum in (2.2). Let $K$ be the unique integer such that $\frac{n_k}{4}\le 2^Kn_k^\beta<\frac{n_k}{2}$. Thus $K = \frac{1-\beta}{\log 2}\log n_k+O(1)$.
Then
\begin{align*}
 \sum_{\frac{n_k}{2}<n_j\le n_k-n_k^\beta} \frac{a_j}{n_k-n_j} & = \sum_{r=1}^K\sum_{2^{r-1}n_k^\beta\le n_k-n_j<2^rn_k^\beta} \frac{a_j}{n_k-n_j} + \sum_{2^Kn_k^\beta\le n_k-n_j<\frac{n_k}{2}} \frac{a_j}{n_k-n_j}\\
 & \le \sum_{r=1}^K\frac{A(n_k-2^{r-1}n_k^\beta) - A(n_k-2^rn_k^\beta)}{2^{r-1}n_k^\beta} + \frac{A(\frac{3n_k}{4}) - A(\frac{n_k}{2})}{2^Kn_k^\beta} \\
 & \ll \sum_{r=1}^K \frac{2^rn_k^\beta}{2^{r-1}n_k^\beta} + 1 \ll K\ll \log n_k.
 \end{align*}
Thus the RHS of (2.2) is
\[ \ll \frac1T\sum_{n_k\leq N}a_kn_k^{1-2\sigma}\log n_k \ll \frac{N^{2-2\sigma}\log N}{T},\]
Taking $T\succ N^{2-2\sigma}\log N$ makes this $o(1)$, as required.

For $\sigma=1$, all the arguments are valid but now the above sum becomes
\[ \sum_{n_k\leq N}\frac{a_k\log n_k}{n_k}\le \log N\sum_{n_k\leq N}\frac{a_k}{n_k}\ll (\log N)^2.\] 
\bo

\noindent
{\bf Proposition 2.2}\nl
{\em Suppose $(1.1)$ holds. Then for each $\delta>0$, we have uniformly for $\sigma\ge\beta+\delta$,}
\[\int_{1}^{T} |f(\sigma+it)-f_N(\sigma+it)|^2\, dt \ll_{\delta} \frac{T^{2}}{N^{2\sigma-2\beta}}+N^{2-2\sigma}.\]
{\em Proof.}\, Write $A(x)=\rho x+R(x)$, where $R(x)\ll x^\beta$. We have for $\sigma>1$,
\begin{align*}
f(s)-f_N(s) & = \sum_{n_j>N}\frac{a_j}{n_j^s} = \int_N^\infty \frac1{x^s}\, dA(x) = \frac{\rho N^{1-s}}{s-1} +  \int_N^\infty \frac1{x^s}\, dR(x)\\
& =\frac{\rho N^{1-s}}{s-1} -\frac{R(N)}{N^s} +s \int_N^\infty \frac{R(x)}{x^{s+1}}\, dx.
\end{align*}
The integral converges absolutely for $\sigma>\beta$ so the above holds in this range by analytic continuation. Thus for $t\in [1,2T]$, we see that 
\[ 	|f(\sigma+it)-f_N(\sigma+it)| \ll \frac{N^{1-\sigma}}{t} + \frac1{N^{\sigma-\beta}} + T\biggl|\int_N^\infty \frac{R(x)}{x^{\sigma+it+1}}\, dx\biggr|.\tag{2.3}\]
Squaring and integrating gives
\[	\int_1^{T} |f(\sigma+it)-f_N(\sigma+it)|^2\, dt 	\ll N^{2-2\sigma} + \frac{T}{N^{2\sigma-2\beta}} + T^2\int_1^{T} \biggl|\int_N^\infty \frac{R(x)}{x^{\sigma+it+1}}\, dx\biggr|^2\, dt.\]
Hence we just need to concentrate on the final term. We shall show that
\[ I(T):=\int_1^{2T} \biggl|\int_N^\infty \frac{R(x)}{x^{\sigma+it+1}}\, dx\biggr|^2\, dt \ll_{\delta} \frac{1}{N^{2\sigma-2\beta}}.\] 
Instead of bounding $I(T)$ directly we use the Fej\'er kernel, and bound
\[
	J(T):=\frac{1}{4T} \int_{-\infty}^{\infty}\Bigl(1-\frac{|t|}{4T}\Bigr)_{+}\Bigl|\int_{N}^{\infty}\frac{R(x)}{x^{\sigma+1+it}}dx\Bigr|^{2}dt,
\]
where $(1-\frac{|t|}{4T})_{+} := \max\{1-\frac{|t|}{4T}, 0\}$.
Note that $I(T) \le 4TJ(T)$. Square out the integral, interchange the order of integration, and use ${\MF\bigl\{\frac{1}{4T}\bigl(1-\frac{|t|}{4T}\bigr)_{+}; u\bigr\}}=s^{2}(2Tu) = \bigl(\frac{\sin(2Tu)}{2Tu}\bigr)^{2}$ with $\MF$ denoting the Fourier transform to see that 
\[
J(T)	\ll \int_{N}^{\infty}\int_{N}^{\infty}\frac{|R(x)R(y)|}{(xy)^{\sigma+1}}s^{2}\bigl(2T\log\tfrac{x}{y}\bigr)d yd x.
\]
By symmetry, it suffices to bound the integral over the domain $N\le y\le x$. We split the integral over $y$ in three parts: $N\le y \le \max\{N, x/2\}$, $\max\{N, x/2\} \le y \le \max\{N, x-x/T\}$, and $\max\{N, x-x/T\} \le y \le x$.

We use $s^{2}(u) \le \min\{1, u^{-2}\}$ and $\log\frac{x}{y} \ge \frac{x-y}{x}$. The contribution of the first range is bounded by
\begin{align*}
	\ll \int_{N}^{\infty}x^{\beta-\sigma-1}\int_{N}^{x/2}\frac{x^{2}y^{\beta-\sigma-1}}{T^{2}(x-y)^{2}}d yd x \ll_{\delta} \int_{N}^{\infty}\frac{x^{\beta-\sigma-1}}{T^{2}N^{\sigma-\beta}}d x \ll_{\delta} \frac{1}{T^{2}N^{2\sigma-2\beta}}.
\end{align*}
The second range gives
\[
	\ll \int_{N}^{\infty}x^{\beta-\sigma-1}\int_{x/2}^{x-x/T}\frac{x^{2}y^{\beta-\sigma-1}}{T^{2}(x-y)^{2}}d yd x \ll \int_{N}^{\infty}\frac{x^{2\beta-2\sigma}}{T^{2}}\cdot\frac{T}{x}d x\ll_{\delta} \frac{1}{TN^{2\sigma-2\beta}}.
\]	
The last range finally yields
\[
	\ll \int_{N}^{\infty}x^{\beta-\sigma-1}\int_{x-x/T}^{x}y^{\beta-\sigma-1}d yd x \ll \int_{N}^{\infty}x^{2\beta-2\sigma-2}\cdot\frac{x}{T}d x \ll_{\delta} \frac{1}{TN^{2\sigma-2\beta}}.
\]
Multiplying by $T$ gives the required bounds for $I(T)$.
\bo

\noindent
{\bf Remark 2}\, It will be useful for later to have a version of Proposition 2.2 for $[T,2T]$, namely:
\[\int_T^{2T} |f(\sigma+it)-f_N(\sigma+it)|^2\, dt \ll_{\delta} \frac{T^{2}}{N^{2\sigma-2\beta}}+\frac{N^{2-2\sigma}}{T}.\]
The proof is identical except that in (2.3), the first term is $\frac{N^{1-\sigma}}{T}$. The discrepancy arises from the fact that $f_N(\sigma+it)$ is large for $t$ small.\nl

For part (b) of Theorem 1.1, we require an extension of Proposition 2.1 for the range $\sigma\le\frac{1+\beta}{2}$. We again employ the Fej\'er kernel, which saves a factor of $\log N$ over the method employed in Proposition 2.1.\nl

\noindent
{\bf Proposition 2.3}\nl
{\em Uniformly for $\sigma\le\frac{1+\beta}{2}$ we have
\[
	\frac{1}{T}\int_{0}^{T}|f_{N}(\sigma+it)|^2 dt \ll \frac{N^{3-2\sigma-\beta}}{T^{2}} + \frac{N^{1+\beta-2\sigma}-1}{1+\beta-2\sigma},
\]
and where the second term of the right hand side is taken to be $\log N$ when $\sigma=\frac{1+\beta}{2}$.}\nl

\noindent
{\em Proof.}\,
Clearly
\[
	\frac{1}{T}\int_{0}^{T}|f_{N}(\sigma+it)|^2 dt \le \frac{1}{T}\int_{-2T}^{2T}\Bigl(1-\frac{|t|}{2T}\Bigr)_{+}|f_{N}(\sigma+it)|^2 dt.
\]
We square out the sum, interchange integration and summation to see that the right hand side equals
\[
	\sum_{n_{k}\le N}\sum_{n_{j}\le N}\frac{a_{k}a_{j}}{(n_{k}n_{j})^{\sigma}}s^{2}\bigl(T\log\tfrac{n_{k}}{n_{j}}\bigr).
\]
By symmetry, it suffices to bound the sum over the range $n_{j}\le n_{k}$. We split this range in three parts: $1\le n_{j}\le n_{k}/2$, $n_{k}/2< n_{j}\le \max\{n_{k}-n_{k}^{\beta}, n_{k}/2\}$, and $\max\{n_{k}-n_{k}^{\beta}, n_{k}/2\} < n_{j}\le n_{k}$. In the first two ranges we bound $s^{2}\bigl(T\log\frac{n_{k}}{n_{j}}\bigr)$ by $\frac{n_{k}^{2}}{T^{2}(n_{k}-n_{j})^{2}}$, while in the last range we bound it by $1$. We obtain
\[
	\sum_{n_{k}\le N}\frac{a_{k}}{n_{k}^{\sigma}}\sum_{n_{j}\le n_{k}/2}\frac{a_{j}}{n_{j}^{\sigma}}\frac{n_{k}^{2}}{T^{2}(n_{k}-n_{j})^{2}} 
	\ll \frac{1}{T^{2}}\sum_{n_{k}\le N}\frac{a_{k}}{n_{k}^{\sigma}}n_{k}^{1-\sigma}\ll \frac{N^{2-2\sigma}}{T^{2}}
\]
for the sum over the first range. The sum over the second range yields
\begin{multline*}
	\sum_{n_{k}\le N}\frac{a_{k}}{n_{k}^{\sigma}}\sum_{n_{k}/2 <n_{j}\le n_{k}-n_{k}^{\beta}}\frac{a_{j}}{n_{j}^{\sigma}}\frac{n_{k}^{2}}{T^{2}(n_{k}-n_{j})^{2}} \\
	\ll\frac{1}{T^{2}}\sum_{n_{k}\le N}\frac{a_{k}}{n_{k}^{\sigma}}n_{k}^{2-\sigma}\sum_{n_{k}/2 <n_{j}\le n_{k}-n_{k}^{\beta}}\frac{a_{j}}{(n_{k}-n_{j})^{2}}
	\ll \frac{1}{T^{2}}\sum_{n_{k}\le N}\frac{a_{k}}{n_{k}^{\sigma}}n_{k}^{2-\sigma-\beta}\ll \frac{N^{3-2\sigma-\beta}}{T^{2}},
\end{multline*}
where we have written 
\[
	\sum_{n_{k}/2 <n_{j}\le n_{k}-n_{k}^{\beta}}\frac{a_{j}}{(n_{k}-n_{j})^{2}} = \int_{n_{k}/2}^{n_{k}-n_{k}^{\beta}}\frac{dA(x)}{(n_{k}-x)^{2}} 
	= \int_{n_{k}/2}^{n_{k}-n_{k}^{\beta}}\frac{d\bigl(\rho x+R(x)\bigr)}{(n_{k}-x)^{2}}
\]
and integrated by parts.
The sum over the last range is bounded as
\[
	\sum_{n_{k}\le N}\frac{a_{k}}{n_{k}^{\sigma}}\sum_{n_{k}-n_{k}^{\beta}<n_{j}\le n_{k}}\frac{a_{j}}{n_{j}^{\sigma}}\ll \sum_{n_{k}\le N}\frac{a_{k}}{n_{k}^{\sigma}}n_{k}^{\beta-\sigma}
	\ll \frac{N^{1+\beta-2\sigma}-1}{1+\beta-2\sigma}.
\]
\bo

\noindent
{\em Proof of Theorem 1.1}. For (a)\, We need only consider $\sigma\le 1$. Combine Propositions 2.1 and 2.2. Take $N$ such that  $N^{2-2\sigma}(\log N)^2\prec T\prec N^{2\sigma-2\beta}$. This is possible since $2-2\sigma<2\sigma-2\beta$ (whenever $\sigma>\frac{1+\beta}{2}$). As such, $\frac1T\int_0^{T} |f-f_N|^2\to 0$. But
\[  
	\frac1T\int_0^{T}|f|^2 = \frac1T\int_0^{T}|f_N|^2+ \frac1T\int_0^{T}|f-f_N|^2 +2\mbox{Re}\biggl(\frac1T\int_0^{T} f_N\overline{(f-f_N)}\biggr).\tag{2.4}
\]
The final term is, in modulus, bounded by $2\bigl(\frac1T\int_0^{T}|f_N|^2 \cdot\frac1T\int_0^{T}|f-f_N|^2\bigr)^{\frac12}$, which tends to 0. By Propositions 2.1 and 2.2, the RHS of (2.4)  converges to the desired limit.

(b)\, Combining Propositions 2.2 and 2.3 with the choice $N=T^{\frac{1}{1-\beta}}$ yields part (b) of Theorem 1.1.
\bo

\noindent
{\bf Generalization to $A(x)=xP(\log x)+O(x^\beta(\log x)^q)$}\nl
Here we consider the same set up but now we replace the main term of (1.1) by $xP(\log x)$ for some polynomial $P$. It is natural to extend the error with an extra $\log$ term too. Assume
\[ A(x) = xP(\log x) +O(x^\beta(\log x)^q)\tag{1.1+}\]
where $P$ is a polynomial of degree $d$ and $q\in\R$. With this condition, $f(s)$ has a pole of order $d+1$ at $s=1$ but retains the analytic continuation for $\sigma>\beta$ as before. Further, we have the bound
\[ A(x)-A(x-y) \ll y(\log x)^d + x^\beta(\log x)^q\] whenever $0\le y< x$, while for $\lambda>-1$, 
\[ \sum_{n_j\le x} a_jn_j^\lambda \sim c_1x^{\lambda+1}(\log x)^d,\quad\mbox{and}\quad
\sum_{n_j\le x} \frac{a_j}{n_j} \sim c_2(\log x)^{d+1},\]
for some $c_1,c_2>0$. Using these, we obtain Propositions 2.1, 2.2, 2.3 and hence Theorem 1.1 with slight modifications, the main one being an extra power of $\log N$ or $\log T$. 
 \nl

\noindent
{\bf Proposition 2.1+}\nl
{\em Suppose $(1.1+)$ holds. Then, for $\sigma\in (\frac{1+\beta}{2},1]$, 
\[ \lim_{T\to\infty}\frac{1}{T}\int_0^T |f_N(\sigma+it)|^2\, dt = \sum_{j=1}^\infty\frac{a_j^2}{n_j^{2\sigma}}\]
as $N,T\to\infty$ such that $T\succ N^{2-2\sigma}(\log N)^{2d+1}$ if $\sigma<1$ and $T\succ(\log N)^{2d+2}$ if $\sigma=1$.}\nl

\noindent
{\bf Proposition 2.2+}\nl
{\em Suppose $(1.1+)$ holds. Then for each $\delta>0$, uniformly for $\sigma\ge\beta+\delta$,}
\[\int_0^T |f(\sigma+it)-f_N(\sigma+it)|^2\, dt \ll_{\delta} \frac{T^2(\log N)^{2q}}{N^{2\sigma-2\beta}}+N^{2-2\sigma}(\log N)^{2d}.\]

\noindent
{\bf Proposition 2.3+}\nl
{\em Uniformly for $\sigma\le\frac{1+\beta}{2}$ we have
\[	\frac{1}{T}\int_{0}^{T}|f_{N}(\sigma+it)|^2 dt \ll \frac{N^{3-2\sigma-\beta}(\log N)^{3d-q}}{T^{2}} + \frac{N^{1+\beta-2\sigma}-1}{1+\beta-2\sigma}(\log N)^{q+d},\]
where the second term of the right hand side is taken to be $(\log N)^{q+d+1}$ when $\sigma=\frac{1+\beta}{2}$.}

\medskip

\noindent
{\bf Theorem 1.1+}\nl
{\em Let $(n_{j})_{j\ge1}$, $(a_{j})_{j\ge1}$ be as above and suppose that (1.1+) holds.
\begin{enumerate}
	\item For\footnote{For $\sigma=1$, we start integrating at $t=1$ to avoid the pole.} $\sigma>\frac{1+\beta}{2}$,  
\[ \lim_{T\to\infty}\frac1T\int_0^T |f(\sigma+it)|^2 \, dt = \sum_{j=1}^\infty\frac{a_j^2}{n_j^{2\sigma}}.	\]
	\item  For $\beta<\sigma<\frac{1+\beta}{2}$, we have
\[\int_{0}^{T}|f(\sigma+it)|^{2}d t \ll T^{\frac{2-2\sigma}{1-\beta}}(\log T)^{\tau(\sigma)},\]
where 
\[
	\tau(\sigma) = 2\frac{(1-\sigma)q+(\sigma-\beta)d}{1-\beta},
\]
while 
\[ 
	\int_{0}^{T}\Bigl|f\Bigl(\frac{1+\beta}{2}+it\Bigr)\Bigr|^{2}d t\ll T(\log T)^{d+q+1}.
\] 
In particular when $q=d$, we have $\tau(\sigma) =2d$.
\end{enumerate}
}

The proofs follow the case where $P(\cdot)$ is constant with minor differences. For example, in the proof of Proposition 2.1+, for the right-most term in (2.2) we now use 
\[ A(n_k-2^{r-1}n_k^\beta(\log n_k)^{q-d}) - A(n_k-2^rn_k^\beta(\log n_k)^{q-d})\ll 2^rn_k^\beta (\log n_k)^q,\]
which leads to a contribution $\ll\frac1T \sum_{n_k\le N} a_kn_k^{1-2\sigma}(\log n_k)^{d+1} \ll \frac{N^{2-2\sigma}(\log N)^{2d+1}}{T}$ (if $\sigma<1$).  

In Proposition 2.2+, we now have 
\[f(s)-f_N(s) = \frac{p_d(\log N,\frac1{s-1})}{N^{s-1}(s-1)} +  \int_N^\infty \frac1{x^s}\, dR(x),\]
where $p_d(x,y)$ is a polynomial in $x$ and $y$ of degree at most $d$. With $t\in [1,2T]$, (2.3) now becomes
\[ 	|f(\sigma+it)-f_N(\sigma+it)| \ll \frac{N^{1-\sigma}(\log N)^d}{t} + \frac{(\log N)^q}{N^{\sigma-\beta}} + T\biggl|\int_N^\infty \frac{R(x)}{x^{\sigma+it+1}}\, dx\biggr|.\]
Squaring and integrating up to $T$ leads to an extra $(\log N)^{2d}$ from the first term, while the bound for $J(T)$ has an extra $(\log N)^{2q}$ (since $R(x)$ has an extra $(\log x)^q$ and the main contribution comes from $N\le x,y\le N^2$).   
Proposition 2.3+ follows a similar pattern.

Finally, for Theorem 1.1+, part (a) follows as before while for (b), take $N=T^{\frac1{1-\beta}}(\log T)^\alpha$ with $\alpha = \frac{q-d}{1-\beta}$.

As an example, $n_j=j$, $a_j=(\log j)^k$ gives $f(s)=(-1)^k\zeta^{(k)}(s)$. We have $\sum_{n\le x} (\log n)^k = xP_k(\log x)+O((\log x)^k)$ for some polynomial $P_k(\cdot)$ of degree $k$, giving $\beta=0$ and $d=q=k$. Theorem 1.1+ says
\[ \lim_{T\to\infty}\frac1T\int_0^T |\zeta^{(k)}(\sigma+it)|^2 \, dt = \sum_{n=1}^\infty\frac{(\log n)^{2k}}{n^{2\sigma}}(=\zeta^{(2k)}(2\sigma)),\]
for $\sigma>\frac12$, while $\int_0^T |\zeta^{(k)}(\sigma+it)|^2 \, dt \ll T^{2-2\sigma}(\log T)^{2k}$ for $\sigma\in (0,\frac12)$ and 
$\int_0^T |\zeta^{(k)}(\frac12+it)|^2 \, dt \ll T(\log T)^{2k+1}$. These are the correct orders of growth, as it was shown by Ingham \cite[Theorem A'' and Theorem A]{Ingham} that 
\begin{align*}
	&\int_{0}^{T}|\zeta^{(k)}(\tfrac{1}{2}+it)|^{2}dt \sim \frac{T}{2k+1}(\log T)^{2k+1}, \\ 
	&\int_{0}^{T}|\zeta^{(k)}(\sigma+it)|^{2}dt \sim \frac{\zeta(2-2\sigma)}{(2\pi)^{1-2\sigma}(2-2\sigma)}T^{2-2\sigma}(\log T)^{2k}, \quad \text{for } \sigma < \frac{1}{2}.
\end{align*}
\nl

\bigskip

\noindent
{\large {\bf 3. Examples}}\nl
In this section we will provide several examples illustrating our results.
\medskip

\noindent
\textbf{Positivity of the coefficients.}
The assumption that the coefficients $a_{j}$ are positive is necessary. Consider for example $n_{j}=j$, $a_{j} = 1 + (-1)^{j}j^{\beta}$. Then $A(x) = x + O(x^{\beta})$, and $f(s) = \zeta(s) + (2^{1+\beta-s}-1)\zeta(s-\beta)$, so that no mean value exists for $\sigma \le \beta+1/2$.
\medskip

\noindent
\textbf{Abscissa of convergence of $\sum_{j} a_{j}^{2}n_{j}^{-2\sigma}$.}
Let $\alpha\in (-\infty, 1)$, and set $n_{j} = j^{\frac{1}{1-\alpha}}$ and $a_{j} = n_{j}^{\alpha} = j^{\frac{\alpha}{1-\alpha}}$.
Then 
\[	
	A(x) =\sum_{j\le x^{1-\alpha}} j^{\frac{\alpha}{1-\alpha}}= (1-\alpha)x + \zeta\Bigl(-\frac{\alpha}{1-\alpha}\Bigr)+ O(x^{\alpha})
\]
so (1.1) holds with $\beta=\alpha$ if $\alpha>0$ and $\beta=0$ otherwise.
The Dirichlet series
\[
	\sum_{j=1}^{\infty}\frac{a_{j}^{2}}{n_{j}^{2\sigma}} = \zeta\Bigl(\frac{2\sigma-2\alpha}{1-\alpha}\Bigr)
\]
has abscissa of convergence $\sigma_{c} = \frac{1+\alpha}{2}$. By adding a sparse sequence $(m_{j})_{j}$ (e.g. $m_{j}=e^{j}$) with coefficients $b_{j}=m_{j}^{\beta}$, it is clear that there exists sequences satisfying
\[
	A(x) = \rho x + O(x^{\beta}), \quad A(x) = \rho x + \Omega(x^{\beta}), 
\]
and for which 
\[
	\sum_{j=1}^{\infty}\frac{a_{j}^{2}}{n_{j}^{2\sigma}}\tag{1.4}
\]
has abscissa of convergence $\sigma_{0}$, for any fixed $\sigma_{0}\in(-\infty, \frac{1+\beta}{2}]$. Setting $n_{j} = \log j$ and $a_{j} = e^{-n_{j}}=1/j$ ($j\ge 2$) with $A(x) = x+\gamma-1+O(e^{-x})$, we can even have that $\sum_{j\ge2}\frac{a_{j}^{2}}{n_{j}^{2\sigma}} = \sum_{j\ge2}(j(\log j)^{\sigma})^{-2}$ has abscissa of convergence $\sigma_{c}=-\infty$.
\medskip

\noindent
\textbf{Existence of the mean value beyond $\sigma=\frac{1+\beta}{2}$.}
We consider the same class of examples as above. If $\alpha < 0$, then $A(x) = (1-\alpha)x+O(1)$, and $f(s) = \zeta\bigl(\frac{s-\alpha}{1-\alpha}\bigr)$ has a mean value if and only if $\sigma>\frac{1+\alpha}{2}$, and $\frac{1+\alpha}{2} < \frac{1}{2}$.

Another example, but with an ordinary Dirichlet series is given by 
\[ 
	f(s) = \sum_{n=1}^\infty \frac{r(n)}{n^s} = 4\zeta(s)\eta(s),
\]
where $r(n)$ counts the number of ways $n$ can be written as a sum of two squares of integers, and where $\eta(s) = 1-3^{-s}+5^{-s}+\cdots$. It is well known that $\sum_{n\le x} r(n)=\pi x +O(x^\theta)$ for some $\theta\in (\frac14,\frac12)$ (indeed it is conjectured that $\theta$ can be taken to be any number greater than $\frac14$). Thus $\beta\ge\frac14$. But it is well-known that the mean value exists for every $\sigma>\frac12$. Here $\frac{1+\beta}{2}\ge \frac58$. Alternatively, we could take $f(s)=\zeta(s)^2$ with $A(x) = \sum_{n\le x}d(n)  = x\log x+(2\gamma-1)x+O(x^\theta)$ for some $\theta\in(\frac14,\frac12)$ (not necessarily the same $\theta$). Theorem 1.1+ says the mean value exists for $\sigma>\frac{1+\theta}{2}$ but, again, it actually exists for every $\sigma>\frac12$.
\medskip

\noindent
\textbf{Sharpness of Theorem 1.1.}
We can employ the same family of basic examples to show that Theorem 1.1 is sharp for every $\beta\in [0,1)$, both with respect to the condition $\sigma>\frac{1+\beta}{2}$, as with respect to the growth bound for $\beta<\sigma\le \frac{1+\beta}{2}$.

Let $n_{j}=j^{\frac{1}{1-\beta}}$ and $a_{j} = n_{j}^{\beta} = j^{\frac{\beta}{1-\beta}}$. Then $f(s) = \zeta\bigl(\frac{s-\beta}{1-\beta}\bigr)$. Using the functional equation, it is straightforward to show that 
\[
	\frac{1}{T}\int_{0}^{T}|\zeta(\sigma+it)|^{2}d t \sim k_\sigma T^{1-2\sigma}, \quad \sigma<1/2,
\]
for some $k_\sigma>0$ (indeed $k_\sigma=\frac{\zeta(2-2\sigma)}{2-2\sigma}(2\pi)^{2\sigma-1}$) so that 
\[
	\frac{1}{T}\int_{0}^{T} |f(\sigma+it)|^{2}d t  \sim\left\{\begin{array}{cl}
		\zeta(\frac{2\sigma-2\beta}{1-\beta}) 	&\text{if } \sigma>\frac{1+\beta}{2};\\
		\log T									&\text{if } \sigma=\frac{1+\beta}{2};\\
		c_\sigma T^{\frac{1+\beta-2\sigma}{1-\beta}}		&\text{if } \sigma<\frac{1+\beta}{2},
	\end{array}\right.
\]
for some $c_\sigma>0$. 
\medskip

An important subclass of Dirichlet series are those with $a_{j}=1$ for every $j$, corresponding to the counting function of the sequence $(n_{j})_{j}$. Also in this subclass, Theorem 1.1 is sharp, as the following example shows. We will further see from this example that convergence of $\sum_{j}a_{j}^{2}n_{j}^{-2\sigma} = \sum_{j}n_{j}^{-2\sigma}$ is not sufficient for the existence of the mean value. It would also be interesting to investigate the sharpness of Theorem 1.1 in the class of Beurling zeta functions.

The example is based on a construction from \cite{BDV}, which in turn was inspired by an example of H.\ Bohr (see e.g.\ the notes to Chapter II.1 of \cite{Tenenbaum}). First we will construct an absolutely continuous non-decreasing function $A_{c}(x)$ for which the bound from Theorem 1.1(b) is sharp, and next we will approximate $A_{c}(x)$ by a function $A(x)$ of the form $A(x) = \sum_{n_{j}\le x}1$ for a certain increasing sequence $(n_{j})_{j}$.

 We fix $\beta \in (0,1)$ and consider a fast growing sequence $(\tau_{k})_{k\ge1}$. We set $B_{k} = \tau_{k}^{\frac{1}{1-\beta}}$ and assume that $\tau_{k}\log B_{k}\in2\pi\Z$ by changing the value of $\tau_{k}$ by a small amount if necessary. We then set
\[
	A_{c}(x) = x-1 + \sum_{k=1}^{\infty}R_{k}(x), \quad x\ge1,
\]
where 
\[
	R_{k}(x) = \left\{\begin{array}{cl}
				\int_{A_{k}}^{x}\cos(\tau_{k}\log u)d u 	&\text{for } A_{k}\le x \le B_{k},\\
				0								&\text{else}.
	\end{array}\right.
\]
Here, $A_{k}$ is chosen in such a way that $R_{k}(B_{k})=0$: we have 
\[
	R_{k}(x) = \biggl[\frac{u\cos(\tau_{k}\log u) - \tau_{k}u\sin(\tau_{k}\log u)}{1+\tau_{k}^{2}}\biggr]^{x}_{A_{k}},
\]
so that we can choose $A_{k}\sim B_{k}/\tau_{k} = B_{k}^{\beta}$ such that $R_{k}(B_{k})=0$. We further assume that $(\tau_{k})_{k}$ increases sufficiently rapidly so that $A_{k+1}>B_{k}$ for $k\ge1$. Clearly, $A_{c}(x)$ is non-decreasing and satisfies $A_{c}(x)=x+O(x^{\beta})$. If we set $f_{c}(s) = \int_{1}^{\infty}x^{-s}d A_{c}(x)$, a small calculation gives
\[
	f_{c}(s) = \frac{1}{s-1} + \frac{1}{2}\sum_{k=1}^{\infty}\Bigl(\frac{B_{k}^{1-s} - A_{k}^{1+i\tau_{k}-s}}{1+i\tau_{k}-s} +\frac{B_{k}^{1-s} - A_{k}^{1-i\tau_{k}-s}}{1-i\tau_{k}-s}\Bigr). 
\]
Here we also used that $\tau_{k}\log B_{k} \in 2\pi\Z$. The singularities at $s = 1 \pm i\tau_{k}$ are removable; we take the value of the corresponding term to be $\log B_{k}/A_{k}$ for $s = 1\pm i\tau_{k}$. If we assume that $\sum_{k}\tau_{k}^{-\eps}<\infty$ for each $\eps>0$, then the above series converges absolutely for $\sigma>\beta$. We also assume now that for each $t>0$, $|\tau_{k}-t| \gg \tau_{k}$ holds for all but at most one value of $k$. Then
\[
	f_{c}(\sigma+i\tau_{k}) = \frac{\tau_{k}^{\frac{1-\sigma}{1-\beta}}}{1-\sigma}\bigl\{1+O\bigl(\tau_{k}^{\sigma-1}\bigr)\bigr\} + O(1), \quad \beta<\sigma<1, \quad f_{c}(1+i\tau_{k})\sim \log\tau_{k}.
\]
In particular, the associated Lindel\"of function $\mu_{f_{c}}(\sigma)$, defined as 
\[
	\mu_{f_{c}}(\sigma) = \inf\{\xi: f_{c}(\sigma+it) \ll |t|^{\xi} \text{ as } |t|\to\infty\},
\]
equals $\mu_{f_{c}}(\sigma) = \bigl(\frac{1-\sigma}{1-\beta}\bigr)_{+}$ for $\sigma>\beta$.

Writing $f(s) = g_{k}(s) + h_{k}(s)$, where $g_{k}(s) = B_{k}^{1-s}/(1+i\tau_{k}-s)$, we have for $\beta<\sigma<1$:
\begin{align*}
	\frac{1}{\tau_{k}}\int_{0}^{\tau_{k}}|f(\sigma+it)|^{2}d t 	&= \frac{1}{\tau_{k}}\int_{0}^{\tau_{k}}|g_{k}(\sigma+it)|^{2}d t 
			+ O\biggl(\frac{1}{\tau_{k}}\int_{0}^{\tau_{k}}|h_{k}(\sigma+it)|^{2}d t\biggr)\\
			&+O\Biggl\{\biggl(\frac{1}{\tau_{k}}\int_{0}^{\tau_{k}}|g_{k}(\sigma+it)|^{2}d t\biggr)^{1/2}\biggl(\frac{1}{\tau_{k}}\int_{0}^{\tau_{k}}|h_{k}(\sigma+it)|^{2}d t\biggr)^{1/2}\Biggr\}.
\end{align*}
It is easily seen that 
\[
	\frac{1}{\tau_{k}}\int_{0}^{\tau_{k}}|g_{k}(\sigma+it)|^{2}d t  \asymp_{\sigma}\tau_{k}^{\frac{1+\beta-2\sigma}{1-\beta}},
\]
while by the rapid increase of $(\tau_{k})_{k}$, 
\begin{align*}
	\frac{1}{\tau_{k}}\int_{0}^{\tau_{k}}|h_{k}(\sigma+it)|^{2}d t 	
		&\ll \tau_{k}^{\beta\frac{2-2\sigma}{1-\beta}-1} + \frac{1}{\tau_{k}}\int_{0}^{2\tau_{k-1}}|h_{k}(\sigma+it)|^{2}d t\\
		&\ll \tau_{k}^{\beta\frac{2-2\sigma}{1-\beta}-1} + 1.
\end{align*}
Hence, when $T=\tau_{k}$, we indeed have
\[
	\frac1T\int_1^T |f_{c}(\sigma+it)|^2 d t \gg_{\sigma} T^{\frac{1+\beta-2\sigma}{1-\beta}}, \quad \text{for } \beta < \sigma < \frac{1+\beta}{2}.
\]

To approximate $A_{c}(x)$ by the counting function of a sequence, we employ a probabilistic method from \cite{BrouckeWeishaupl}.\nl

\noindent
{\bf Theorem} {\cite[Theorem 2.3]{BrouckeWeishaupl}}\nl
{\em Let $M: [1,\infty) \to [0,\infty)$ be a continuous and non-decreasing function satisfying $M(1)=0$ and $\int_{1}^{x}\sqrt{M(u)}d u/u \ll \sqrt{M(x)}$. Set $x_{j} = \inf\{x: M(x)=j\}$. select $n_{j}$ randomly and independently from the interval $(x_{j-1}, x_{j}]$ according to the probability distribution $d M(x)\big|_{(x_{j-1}, x_{j}]}$. Then with probability $1$ we have for each $x\ge1$ and $t\in\R$:}
\[	\biggl| \sum_{n_{j}\le x}n_{j}^{-i t} - \int_{1}^{x}u^{-i t}d M(u) \biggr| \ll \sqrt{M(x)}\bigl(\sqrt{\log(x+1)} + \sqrt{\log(\,|t|+1)}\bigr).\tag{3.1}\]
We apply the theorem with $M=A_{c}$ to find a sequence $(n_{j})_{j}$ satisfying the above bound. By construction of the random variables, 
\[
	A(x) := \sum_{n_{j}\le x}1 = A_{c}(x) + O(1) = x + O(x^{\beta}).
\]
We set
\[
	f(s) = \sum_{j=1}^{\infty}n_{j}^{-s}, \quad \sigma>1, \quad S_{t}(x) = \sum_{n_{j}\le x}n_{j}^{-it}, \quad I_{t}(x) = \int_{1}^{x}u^{-it}d A_{c}(u).
\]
Clearly $f(s) - 1/(s-1)$ has analytic continuation to $\sigma>\beta$, and $f(s)-f_{c}(s)$ has analytic continuation to $\sigma>0$. Employing the bound (3.1) and integrating by parts, we have for $\sigma>1/2$ that
\begin{align*}
	f(s) -f_{c}(s) 	&= \int_{1}^{\infty}x^{-\sigma}d\,\bigl(S_{t}(x) - I_{t}(x)\bigr) = \sigma\int_{1}^{\infty}x^{-\sigma-1}\bigl(S_{t}(x) - I_{t}(x)\bigr)d x \\
				&\ll \int_{1}^{\infty}x^{-\sigma-1/2}\bigl(\sqrt{\log x}+\sqrt{\log |t|}\bigr)d x = \frac{1}{(\sigma-1/2)^{3/2}} + \frac{\sqrt{\log |t|}}{\sigma-1/2}.
\end{align*}
For $0< \sigma\le 1/2$, we will interpolate the bounds resulting from (3.1) and $A(x) - A_{c}(x) \ll 1$. Let $0<\sigma\le 1/2$ and write
\[
	f(s) - f_{c}(s) = \int_{1}^{T}x^{-\sigma}d\,\bigl(S_{t}(x) - I_{t}(x)\bigr) + \int_{T}^{\infty}x^{-s}d\,\bigl(A(x) - A_{c}(x)\bigr).
\]
The first term is $\ll \frac{T^{1/2-\sigma}-1}{1/2-\sigma}\sqrt{\log T}$ (which we take to be $(\log T)^{3/2}$ for $\sigma=1/2$), while for the second term we have
\[
	\int_{T}^{\infty}x^{-s}d\,\bigl(A(x)-A_{c}(x)\bigr) = O(T^{-\sigma})+s\int_{T}^{\infty}\frac{A(x)-A_{c}(x)}{x^{\sigma+it+1}}dx.
\]
Hence,
\begin{align*}
	\int_{0}^{T}|f(\sigma+it)-f_{c}(\sigma+it)|^{2}dt 	&\ll T^{2-2\sigma}(\log T)^{3} + T^{2}\int_{0}^{T}\biggl|\int_{T}^{\infty}\frac{A(x)-A_{c}(x)}{x^{\sigma+it+1}}dx\biggr|^{2}dt\\
			&\ll T^{2-2\sigma}(\log T)^{3} + \frac{T^{2-2\sigma}}{\sigma}+\frac{T^{1-2\sigma}}{\sigma^{2}}.
\end{align*}
Indeed, the last term of the first line can be estimated in the same way as the estimation of $I(T)$ from the proof of Proposition 2.2. The obtained estimates now imply that also 
\[
	\frac{1}{\tau_{k}}\int_{0}^{\tau_{k}}|f(\sigma+it)|^{2}d t \gg_{\sigma} \tau_{k}^{\frac{1+\beta-2\sigma}{1-\beta}}, \quad \text{ for } \beta <\sigma < \frac{1+\beta}{2}.
\]

Finally we note that this example also shows that convergence of $\sum_{j}a_{j}^{2}n_{j}^{-2\sigma}$ is not sufficient for the existence of the mean value: here $\sum_{j} a_{j}^{2}n_{j}^{-2\sigma}$ converges if and only if $\sigma>\frac{1}{2}$, while the mean value exists only if $\sigma > \frac{1+\beta}{2}$.
%
%
\bigskip

\noindent
{\large {\bf 4. Behavior on the line $\sigma=(1+\beta)/2$}}\nl
Theorem 1.1(b) yields the upper bound $\int_{0}^{T}|f(\frac{1+\beta}{2}+it)|^{2}dt \ll T\log T$ on the ``critical line'' $\sigma=\frac{1+\beta}{2}$. In the case of the Riemann zeta function, we actually have the asymptotic $\int_{0}^{T}|\zeta(\frac{1}{2}+it)|^{2}dt \sim T\log T$. We first show, assuming some separation of the sequence $(n_{j})_{j}$, a similar asymptotic on the line $\sigma=\frac{1+\beta}{2}$. For this, we use the following mean value theorem for ``general Dirichlet polynomials'' due to Montgomery and Vaughan, \cite[Corollary 2]{MV74}.\nl

\noindent
{\bf Theorem} [Montgomery, Vaughan]\nl
{\em Suppose that $\lambda_{1}, \lambda_{2},\dotsc, \lambda_{J}$ are distinct real numbers, and that $a_{j}$, $j=1,\dotsc, J$ are complex numbers. If $\delta_{j} = \min_{k\neq j}|\lambda_{k}-\lambda_{j}|$, then 
\[
	\int_{0}^{T}\biggl|\sum_{j=1}^{J}a_{j}\exp(i\lambda_{j}t)\biggr|^{2}dt = \sum_{j=1}^{J}|a_{j}|^{2}(T+ 3\pi\theta\delta_{j}^{-1}), 
\]
where $\theta$ is a real number depending on the various parameters satisfying $|\theta|\le 1$}.\nl

\noindent
{\bf Theorem 4.1}\nl
{\em Suppose that $n_{j+1}-n_{j} \gg \max\{a_{j}, a_{j+1}\}$. Then}
\[
	\frac{1}{T}\int_{0}^{T}\Bigl|f\Bigl(\frac{1+\beta}{2}+it\Bigr)\Bigr|^{2}d t = \sum_{n_{j}\le T^{\frac{1}{1-\beta}}}\frac{a_{j}^{2}}{n_{j}^{1+\beta}} + O\Biggl(\biggl(\sum_{n_{j}\le T^{\frac{1}{1-\beta}}}\frac{a_{j}^{2}}{n_{j}^{1+\beta}}\biggr)^{1/2}\Biggr).
\]
Note that in case of convergence of the series on the right, this reduces to $\frac{1}{T}\int_{0}^{T}|\dotso|^{2}=O(1)$. In any case, as $a_{j}\ll n_{j}^{\beta}$, the series is $\ll \log T$.

\noindent
{\em Proof.}
Setting $N=T^{\frac{1}{1-\beta}}$, we have by Proposition 2.2
\[
	\frac{1}{T}\int_{0}^{T}\Bigl|(f-f_{N})\Bigl(\frac{1+\beta}{2}+it\Bigr)\Bigr|^{2}d t \ll 1.
\]
Applying the Montgomery--Vaughan Theorem, we get that
\[
	\frac{1}{T}\int_{0}^{T}\Bigl|f_{N}\Bigl(\frac{1+\beta}{2}+it\Bigr)\Bigr|^{2}d t = \sum_{n_{j}\le N}\bigl(1+O(T\delta_{j})^{-1}\bigr)\frac{a_{j}^{2}}{n_{j}^{1+\beta}},
\]
where $\delta_{j}:= \min_{k\neq j}|\log n_{j}-\log n_{k}|$, which is $\gg a_{j}/n_{j}$ by assumption. Hence the big-$O$ term above is 
\[
	\ll\frac{1}{T}\sum_{n_{j}\le N}\frac{a_{j}}{n_{j}^{\beta}} \ll \frac{N^{1-\beta}}{T}\ll 1.
\]
The result now follows from Cauchy--Schwartz: 
\[
	\int_{0}^{T} |f|^{2} = \int_{0}^{T}|f_{N}|^{2} + \int_{0}^{T}|f-f_{N}|^{2} + O\Biggl(\biggl(\int_{0}^{T}|f_{N}|^{2}\int_{0}^{T} |f-f_{N}|^{2}\biggr)^{1/2}\Biggr).
\]
\bo
Let us highlight the particular case $\beta=0$, $a_{j}=1$. If one is merely interested in the asymptotic $\sim \rho \log T$ of the mean square on the line $\sigma=1/2$, one can ask for a slightly more general separation condition.\nl

\noindent
{\bf Corollary 4.2}\nl
{\em Let $g: [0,\infty) \to [0,1]$ be a non-increasing function satisfying $1/g(x) = o(\log^{2}x)$. Let $(n_{j})_{j\ge1}$ be an increasing sequence of real numbers, $n_{1}\ge1$, satisfying $n_{j+1}-n_{j}\gg g(n_{j})$ ($j\ge1$) and $N(x) = \sum_{n_{j}\le x}1 = \rho x + O(1)$ for some $\rho>0$. Set
\[
	f(s) = \sum_{j\ge1}\frac{1}{n_{j}^{s}},
\]
initially defined for $\Re s>1$, and by analytic continuation to $\Re s>0$ with the exception of a simple pole at $s=1$ with residue $\rho$. Then}
\[
	\frac{1}{T}\int_{0}^{T}|f(\sigma+it)|^{2}d t \sim 
	\left\{\begin{array}{cl}
		f(2\sigma)	&\text{if } \sigma>\frac{1}{2};\\
		\rho\log T	&\text{if } \sigma=\frac{1}{2}.
	\end{array}\right.
\]

\noindent
{\em Proof.}
Following the above proof using the Montgomery--Vaughan Theorem, we get, assuming without loss of generality that $g(n_{j+1}) \gg g(n_{j})$,
\[
	\frac{1}{T}\int_{0}^{T}\Bigl|f_{N}\Bigl(\frac{1}{2}+it\Bigr)\Bigr|^{2}d t = \sum_{n_{j}\le N}\Bigl(1+O\Bigl(\frac{n_{j}}{Tg(n_{j})}\Bigr)\Bigr)\frac{1}{n_{j}}.
\]
The contribution from the $O$-term is bounded by
\[
	\ll \sum_{n_{j}\le N}\frac{1}{Tg(n_{j})} \ll \frac{N}{Tg(N)}.
\]
Hence we require that $N$ is such that $N/g(N) \prec T\log T$. When invoking Proposition 2.2 with $\beta=0$ and $\sigma=\frac12$, we also require that $T/\log T \prec N \prec T\log T$. Taking $N=T\sqrt{g(T)}$, this is satisfied, while $N\le T$ implies $N/g(N)=T\sqrt{g(T)}/g(N)\le T/\sqrt{g(T)}=o(T\log T)$.
\bo

Without sufficient separation of the $n_{j}$, the contribution of the off-diagonal terms can become significant and the result may fail.\nl

\noindent
{\bf Proposition 4.3}\nl
{\em Let $\delta>0$ be arbitrary. Then there exists an increasing sequence of real numbers $(n_{j})_{j\ge 1}$ satisfying $N(x) = \sum_{n_{j}\le x}1 = \rho x + O(1)$ for some $\rho>0$, $n_{j+1}-n_{j} \gg n_{j}^{-\delta}$, and for which
\[
	\frac{1}{T}\int_{0}^{T}\Bigl|f\Bigl(\frac{1}{2}+it\Bigr)\Bigr|^{2}d t \not\sim \rho \log T.
\]} 

\noindent
{\em Proof.}
Let $L$ be a large fixed integer and set $m_{k,l} = k + lk^{-\delta}$ for $k\ge k_{0}$, $l=0,1,\dotsc, L-1$. Here $k_{0}$ is such that $k+Lk^{-\delta} < k+1$ if $k\ge k_{0}$. We then set $n_{j} = m_{k,l}$ if $j=(k-k_{0})L + l +1$; clearly $N(x)  = Lx+O(1)$. We have 
\begin{align*}
	\frac{1}{T}\int_{0}^{T}\Bigl|f_{T}\Bigl(\frac{1}{2}+it\Bigr)\Bigr|^{2}d t 
		&\ge \frac{1}{2T}\int_{-T}^{T}\Bigl(1-\frac{|t|}{T}\Bigr)_{+}\Bigl|f_{T}\Bigl(\frac{1}{2}+it\Bigr)\Bigr|^{2}d t = \frac{1}{2}\sum_{n_{j}, n_{r}\le T}\frac{s^{2}\bigl(\tfrac{T}{2}\log\tfrac{n_{j}}{n_{r}}\bigr)}{(n_{r}n_{s})^{1/2}}\\
		&\ge \frac{1}{2}\sum_{n_{j}\le T}\frac{1}{n_{j}} + \sum_{k\le T/L-1}\sum_{l_{1}<l_{2}}\frac{s^{2}\bigl(\tfrac{T}{2}\log\tfrac{m_{k,l_{1}}}{m_{k,l_{2}}}\bigr)}{k+1}.\\
\end{align*}
Now $\log\frac{m_{k,l_{1}}}{m_{k,l_{2}}} \sim (l_{1}-l_{2})k^{-1-\delta}$, so there is a constant $c$ such $s^{2}\bigl(\tfrac{T}{2}\log\tfrac{m_{k,l_{1}}}{m_{k,l_{2}}}\bigr) \ge 1/2$ whenever $k\ge cT^{\frac{1}{1+\delta}}$.
We get 
\begin{align*}
	\frac{1}{T}\int_{0}^{T}\Bigl|f_{T}\Bigl(\frac{1}{2}+it\Bigr)\Bigr|^{2}d t  	
			&\ge \frac{1}{2}\sum_{n_{j}\le T}\frac{1}{n_{j}} + \sum_{cT^{\frac{1}{1+\delta}}\le k \le T/L-1} \sum_{l_{1}<l_{2}}\frac{1}{2(k+1)}\\
			&\ge \frac{L}{2}\log T + \frac{\delta}{1+\delta}\frac{L(L-1)}{4}\log T + o(\log T).
\end{align*}
The claim follows upon choosing $L$ sufficiently large in terms of $\delta$, and applying Proposition 2.2.
\bo
It is an interesting question how small a ``gap function'' $g$ can be so that $n_{j+1}-n_{j}\gg g(n_{j})$ implies that the mean square of $f$ on $\sigma =1/2$ is asymptotic to $\rho \log T$. We have shown that $g(x) \succ 1/\log^{2}x$ suffices but that $g(x) \gg x^{-\delta}$ is not sufficient, for each $\delta>0$. 
%

\medskip

Theorem 4.1 can also be used to show that in general
\[ \frac1{T\log T} \int_0^T\Bigl|f\Bigl(\frac{1+\beta}{2}+it\Bigr)\Bigr|^{2}dt\]
does not converge as $T\to\infty$, even for ordinary Dirichlet series. We shall consider the case $\beta=0$ (a simple adjustment can be done for other values of $\beta$). Let $n_j=j$ and $a_j=1+\eta_j$, where 
$\eta_j$ is chosen so that $\sum_{j\le x} \eta_j  =0$ or 1. Indeed, let $M_n=2^{2^n}$ for $n\in\N_0$ and define
\[ \eta_n = \left\{ \begin{array}{cl} (-1)^n & \mbox{ for $M_{2N}\le n< M_{2N+1}$}\\ 0 & \mbox{ for $M_{2N+1}\le n< M_{2N+2}$}\end{array}\right. .\]
Thus $A(x) = x+O(1)$ and 
\[ \sum_{n_{j}\le T^{\frac{1}{1-\beta}}}\frac{a_{j}^{2}}{n_{j}^{1+\beta}} = \sum_{j\le T} \frac{1+2\eta_j+\eta_j^2}{j} = \log T + \sum_{j\le T} \frac{\eta_j^2}{j} +O(1).\]
But the sum on the right is not asymptotic to $c\log T$ since
\[ \sum_{j=M_{2N}}^{M_{2N+1}} \frac{\eta_j^2}{j} = \sum_{j=M_{2N}}^{M_{2N+1}} \frac{1}{j} \sim \frac{1}{2}\log M_{2N+1}\quad \mbox{ while }\quad \sum_{j=M_{2N+1}}^{M_{2N+2}-1} \frac{\eta_j^2}{j} =0.\]
For the general case, take $n_j=j^{\frac1{1-\beta}}$ and $a_j=n_j^\beta(1+\eta_j)$, with $\eta_j$ as above. Furthermore:\nl

\noindent
{\bf Proposition 4.4}\nl
{\em For the above example, 
\[ \frac1{T^{2-2\sigma}} \int_0^T |f(\sigma+it)|^{2}dt\]
does not converge as $T\to\infty$, for any $0<\sigma<\frac12$.

\medskip
\noindent
Proof.}\, Write $f$ for $f(\sigma+it)$ and $f_N$ for $f_N(\sigma+it)$. 

Fix $\sigma\in (0,\frac12)$ and let $0<a<1<A$. By Proposition 2.3, we have
\[ \int_T^{2T} |f_{aT}|^2 \ll a^{1-2\sigma}T^{2-2\sigma},\]
while by the Remark following Proposition 2.2,
\[ \int_T^{2T} |f-f_{AT}|^2 \ll \frac{T^{2-2\sigma}}{A^{2\sigma}},\]
both uniform in $a,A$, for $T\ge T_0(A)$, some $T_0(A)$. Hence, given $\eps>0$, there exist $a_0,A_0>0$ such that 
$\int_T^{2T} |f_{aT}|^2 +|f-f_{AT}|^2<\eps T^{2-2\sigma}$ for $a\le a_{0}$, $A\ge A_{0}$.

Let $\Sigma = f_{AT} - f_{aT} = \sum_{aT<n\le AT} \frac{a_n}{n^{\sigma+it}}$. Then $\int_T^{2T} |f-\Sigma|^2<2\eps T^{2-2\sigma}$. It follows that (using $|z+w|^2\le (1+\sqrt{\eps})|z|^2 + (1+\frac1{\sqrt{\eps}})|w|^2$) 
\[ \int_T^{2T} |\Sigma|^2\le (1+\sqrt{\eps})\int_T^{2T} |f|^2 + 4\sqrt{\eps}T^{2-2\sigma}.\]
Similarly vice versa, with $\Sigma$ and $f$ swapped. 

Now suppose, $\int_T^{2T} |f|^2 \sim cT^{2-2\sigma}$ for some $c>0$. Then, for all $\delta>0$, there exist $a^\prime,A^\prime>0$ such that 
\[ \frac1{T^{2-2\sigma}}\int_T^{2T} \biggl|\sum_{aT<n\le AT}\frac{a_n}{n^{\sigma+it}}\biggr|^2 \in (c-\delta,c+\delta)\tag{4.1}\]
for $T$ large enough and any $a\le a^\prime$ and $A\ge A^\prime$. 

By the same method and using the fact that $\int_T^{2T}|\zeta(\sigma+it)|^2\, dt \sim dT^{2-2\sigma}$ for some $d>0$, we have
\[ \frac1{T^{2-2\sigma}}\int_T^{2T} \biggl|\sum_{aT<n\le AT}\frac1{n^{\sigma+it}}\biggr|^2 \in (d-\delta,d+\delta)\tag{4.2}\]
for $T$ sufficiently large whenever $a\le a^\prime$ and $A\ge A^\prime$.

But (4.1) and (4.2) are in contradiction. Choose $a=a^\prime$ and $A=2A^\prime$ and take $T$ of the form
\[ T=2^{2^{k+\lambda}}\]
with $k\in\N_0$ and $\lambda\in [0,2)$. If $\lambda\in (1,2)$ then, for $k$ large enough, $[aT,AT]\subset [M_{2k+1},M_{2k+2})$ and so $a_n=1$ for all $n\in [aT,AT]$. Thus $\Sigma = \sum_{aT<n\le AT}n^{-\sigma-it}$ and $|c-d|<2\delta$ must follow.

On the other hand, if $0<\lambda<1$, then for $T$ large enough, $[aT,AT]\subset [M_{2k},M_{2k+1})$ and so $a_n=1+(-1)^n$ for all $n\in [aT,AT]$. Thus $a_n=2$ for $n$ even and zero otherwise in this interval, so that
\[\Sigma = \sum_{aT<2n\le AT}\frac{2}{(2n)^{\sigma+it}} = 2^{1-\sigma}\sum_{\frac{aT}{2}<n\le \frac{AT}{2}}\frac1{n^{\sigma+it}}\]
and $|c-2^{2-2\sigma}d|<2\delta$ must also hold. For $\delta$ small enough, this is incompatible with $|c-d|<2\delta$.
\bo


\bigskip

\noindent
{\large {\bf 5. The number of zeros of general Dirichlet series}}\nl
In 1913, Bohr and Landau \cite{BohrLandau} observed that the boundedness of the mean square of the Riemann zeta function $\zeta(s)$ for $\sigma>1/2$ implies that its number of zeros up to height $T$ in half-planes to the right of the critical line $\sigma=1/2$ is bounded by a constant times $T$. More precisely, setting
\[
	N(\sigma, T) = \#\{s'=\sigma'+it': \zeta(s')=0,\, \sigma'\ge\sigma,\, |t'|\le T\},
\]
they showed that $N(\sigma,T) \ll_{\sigma} T$ for $\sigma>1/2$. 
Armed with our new mean value theorem, we can now show a similar statement about the zeros of general Dirichlet series with positive coefficients. Incidently, this answers a question posed in \cite{BrouckeDebruyne}. In that paper, a zero-density estimate was obtained for the number of zeros of Beurling zeta functions, that is, general Dirichlet series of the form 
\[
	\zeta_{\MP}(s) = \sum_{j=1}^{\infty}\frac{a_{j}}{n_{j}^{s}},
\]
where $(n_{j})_{j}$ is the sequence of generalized integers from a Beurling generalized number system $(\MP, \MN)$ with corresponding multiplicities $a_{j}$. The obtained zero-density estimate holds in the range $\sigma > \frac{2+\beta}{3}$, in case that $\beta$ is such that $\sum_{n_{j}\le x}a_{j} = \rho x + O(x^{\beta})$ holds for some $\rho>0$. It was asked whether anything better than the ``trivial'' estimate $N(\sigma,T) \ll T\log T$ could be obtained in the range $\frac{1+\beta}{2} < \sigma \le \frac{2+\beta}{3}$. The following theorem provides a positive answer.\nl

\noindent
{\bf Theorem 5.1}\nl
{\em Let $f(s)  = \sum_{j}a_{j}n_{j}^{-s}$ be a general Dirichlet series with positive coefficients satisfying $(1.1)$. Given $\sigma>\beta$, $T>0$, let $N(\sigma, T)$ be the number of zeros $s'=\sigma'+it'$ of $f(s)$ satisfying $\sigma'\ge\sigma$ and $|t'|\le T$. Then for $\sigma> \frac{1+\beta}{2}$, 
\[
	N(\sigma, T) \ll_{\sigma} T.
\]
Proof.}\, The proof is very similar to the proof in the case of ordinary Dirichlet series having finite order and a mean value, see e.g.\ \cite[\S 9.623]{Titchmarsh}. We give an overview here, but refer to \cite{Titchmarsh} for more details.

Notice first that $f(s) = a_{1}n_{1}^{-s} + O(n_{2}^{-\sigma})$ for $\sigma$ sufficiently large. Fix $\sigma_{1}$ so that $C_{1} < |f(\sigma_{1}+ it)| < C_{2}$ for certain positive constants $C_{1}$, $C_{2}$, and so that $f(s)$ has no zeros in the half-plane $\sigma\ge\sigma_{1}$. Assume that $f$ has no zeros on the rectangle with vertices $\sigma\pm iT$, $\sigma_{1}\pm iT$, by slightly decreasing the value of $\sigma$ and increasing the value of $T$ if necessary. By a formula of Littlewood, 
\[
	2\pi\int_{\sigma}^{\sigma_{1}}\bigl(N(\alpha,T)-\chi_{[1,\infty)}(\alpha)\bigr)d\alpha = \int_{-T}^{T}\log|f(\sigma+it)|dt - \int_{-T}^{T}\log|f(\sigma_{1}+it)|dt + 2\int_{\sigma}^{\sigma_{1}}\arg f(\alpha+iT)d\alpha.
\]
Here, $\chi_{[1,\infty)}$ is the characteristic function of $[1, \infty)$, representing the contribution of the pole at $s=1$, and the branch of $\arg f(s)$ is chosen which is zero when $s=\sigma>1$ is real (when also $f(s)$ is real), and which is continuous along the broken line $[\sigma_{1},\sigma_{1}+iT]\cup[\sigma+iT,\sigma_{1}+iT]$.

The second integral is plainly $O(T)$, as $\log|f(\sigma_{1}+it)| \ll 1$. For the third integral, we observe that $\arg f(\sigma_{1}+iT) \ll T$. 
The variation of $\arg f(s)$ along the line $[\sigma+iT, \sigma_{1}+iT]$ is bounded by a constant times the number of zeros of $\Re f$ along this line. Using Jensen's theorem and the fact that $f$ has finite order, one sees that this variation is $O(\log T)$, so that also the third integral is $O(T)$. Finally for the first integral, we use Jensen's inequality and the existence of the mean to see that
\[
	\int_{-T}^{T}\log|f(\sigma+it)|dt = T\int_{-T}^{T}\log|f(\sigma+it)|^{2}\frac{dt}{2T} \le T\log\biggl(\int_{-T}^{T}|f(\sigma+it)|^{2}\frac{dt}{2T}\biggr) \ll T.
\]
\bo
%

Of course, as evident from the Riemann zeta function, the condition $\sigma> \frac{1+\beta}{2}$ is sharp.

\bigskip

\noindent
{\bf Appendix: (A)  General Mellin--Stieltjes transforms}\nl
Theorem 1.1 can be extended from Dirichlet series to more general Mellin--Stieltjes transforms. Let $A: [0,\infty) \to [0,\infty)$ be a non-decreasing and right-continuous function with $A(x)=0$ for $0\le x<1$. Suppose $A(x) = \rho x + O(x^{\beta})$ for some $\rho>0$ and $0\le\beta<1$. Set 
\[
	f(s) = \int_{0}^{\infty}x^{-s}d A(x).
\]
Then Theorem 1.1 holds, upon replacing $\sum_{j}\frac{a_{j}^{2}}{n_{j}^{2\sigma}}$ by
\[
	\iint_{x=y}(xy)^{-\sigma}d A(x)d A(y).
\]

The above integral actually reduces to a series. To see this, we observe that the product of two positive measures vanishes on the diagonal as soon as one of them is continuous: let $F$, $G$ be non-decreasing right-continuous functions, and suppose $F$ is continuous. Then for all $x\ge0$, $d F\times d G(\{(u,u): 0\le u\le x\}) = 0$. Indeed, let $\eps>0$ be arbitrary. As $F$ is uniformly continuous on $[0,x]$, we may take a $\delta>0$ sufficiently small such that $|u-v|\le\delta \implies |F(u)-F(v)| < \eps$. Covering the diagonal with small squares, and assuming that $x/\delta$ is an integer (which we may by making $\delta$ smaller), we have
\begin{align*}
	d F\times d G	&(\{(u,u): 0\le u\le x\})\\
					&\le (F(\delta)-F(0))(G(\delta)-G(0)) + \dotsb + (F(x)-F(x-\delta))(G(x)-G(x-\delta))\\
					&\le \eps\bigl(G(\delta)-G(0)+\dotsb+G(x)-G(x-\delta)\bigr)=\eps(G(x)-G(0)).
\end{align*}
But $\eps>0$ was arbitrary, so the claim follows.

Let $(n_{j})_{j\ge1}$ be the increasing sequence of points of discontinuity of $A$. Then 
\[
	A(x) = \sum_{n_{j}\le x}a_{j} + A_{c}(x)
\]
for certain positive real numbers $a_{j}$ and a continuous non-decreasing function $A_{c}$. By the observation, we obtain
\[
	\iint_{x=y}(xy)^{-\sigma}d A(x)d A(y) = \sum_{j=1}^{\infty}\frac{a_{j}^{2}}{n_{j}^{2\sigma}}.
\]

\medskip

\noindent
{\bf (B) Bessel's inequality}\nl
Here we show (1.3) (see also e.g.\ \cite[p. 11]{CI}).  
By absolute convergence, we have for $\sigma>1$
\[ \frac1{2T}\int_{-T}^T f(\sigma+it)n_j^{it}\, dt \to a_jn_j^{-\sigma}\tag{B1}\]
as $T\to\infty$. Equivalently, 
 \[ \frac1{2iT}\int_{\sigma-iT}^{\sigma+iT} f(s)n_j^s\, ds \to a_j.\tag{B2}\]
But the analytic continuation of $f$ implied by (1.1) also gives $f(\sigma+it)=o(t)$ uniformly for $\sigma\ge\beta+\delta$ for any $\delta>0$. It follows that (B2) and hence (B1) also holds for every $\sigma>\beta$. For any fixed $N$, we therefore have with $f_{N}(s) = \sum_{n_{j}\le N}a_{j}n_{j}^{-s}$ that
\[ \frac1{2T}\int_{-T}^T f\overline{f_N} = \sum_{n_j\le N}\frac{a_j}{n_j^\sigma}\frac1{2T}\int_{-T}^T f(\sigma+it)n_j^{it}\, dt \to \sum_{n_j\le N}\frac{a_j^2}{n_j^{2\sigma}},\]
as $T\to\infty$. Thus
\begin{align*} 
0 &\le \frac1{2T}\int_{-T}^T |f-f_N|^2 = \frac1{2T}\int_{-T}^T |f|^2 + \frac1{2T}\int_{-T}^T |f_N|^2 - 2 \mbox{Re}\biggl( \frac1{2T}\int_{-T}^T f\overline{f_N}\biggr)\\
& =\frac1{2T}\int_{-T}^T |f|^2 - \sum_{n_j\le N}\frac{a_j^2}{n_j^{2\sigma}}+o(1).
\end{align*}
As $N$ is arbitrary, it follows that 
\[ \liminf_{T\to\infty}\frac{1}{T}\int_0^T |f(\sigma+it)|^2\, dt \ge \sum_{j=1}^{\infty}\frac{a_{j}^{2}}{n_{j}^{2\sigma}}.\]

\noindent
{\bf Funding}\nl
This work was supported by the Research Foundation -- Flanders [12ZZH23N to F.B.].

\noindent
Frederik Broucke, Department of Mathematics: Analysis, Logic and Discrete Mathematics, Ghent University, Krijgslaan 281, 9000 Gent, Belgium.\\
Email address: fabrouck.broucke@ugent.be

\medskip
\noindent
Titus W. Hilberdink, Nanjing University of Information Science and Technology (Reading Academy), 219 Ningliu Road, Nanjing, China.\\
Email address: t.w.hilberdink@nuist.edu.cn, t.hilberdink@reading.ac.uk

\end{document}